\documentclass[12pt]{amsart}
\usepackage{a4wide,enumerate,xcolor}
\usepackage{amsmath,graphicx,comment}
\usepackage{esint} 
\allowdisplaybreaks

\let\pa\partial
\let\na\nabla

\newcommand{\R}{{\mathbb R}}
\newcommand{\diver}{\operatorname{div}}

\newtheorem{theorem}{Theorem}
\newtheorem{lemma}[theorem]{Lemma}

\newtheorem{remark}[theorem]{Remark}
\newtheorem{corollary}[theorem]{Corollary}
\newtheorem{definition}{Definition}

\begin{document}

\title[Stability and uniqueness of triangular cross-diffusion systems]{Stability and uniqueness of bounded weak solutions to triangular degenerate cross-diffusion systems} 

\author[X. Chen]{Xiuqing Chen}
\address{School of Mathematics (Zhuhai), Sun Yat-sen University, Zhuhai 519082, Guangdong Province, China}
\email{chenxiuqing@mail.sysu.edu.cn}
    
\author[B. Du]{Bang Du}
\address{School of Mathematics (Zhuhai), Sun Yat-sen University, Zhuhai 519082, Guangdong Province, China}
\email{dubang5@mail.sysu.edu.cn}

\author[A. J\"ungel]{Ansgar J\"ungel}
\address{Institute of Analysis and Scientific Computing, TU Wien, Wiedner Hauptstra\ss e 8--10, 1040 Wien, Austria}
\email{juengel@tuwien.ac.at} 

\date{\today}

\thanks{The first and second authors acknowledge support from the National Natural Science Foundation of China (NSFC),  grant 12471206. The last author acknowledges partial support from the Austrian Science Fund (FWF), grant 10.55776/F65, and from the Austrian Federal Ministry for Women, Science and Research and implemented by \"OAD, project MultHeFlo. This work has received funding from the European Research Council (ERC) under the European Union's Horizon 2020 research and innovation programme, ERC Advanced Grant NEUROMORPH, no.~101018153. For open-access purposes, the authors have applied a CC BY public copyright license to any author-accepted manuscript version arising from this submission.} 

\begin{abstract}
The continuous dependence on the initial data and consequently the uniqueness of bounded weak solutions to a class of triangular reaction-cross-diffusion equations is shown. The class includes two-species doubly degenerate equations for nutrient taxis models describing the response of bacteria to nutrient conditions. The key difficulty is the lack of a gradient bound for the difference of the first component of the solution, due to the degeneracy. This issue is overcome by assuming a nonstandard Lipschitz-type condition, applying the $H^{-1}$ method, and carefully combining various estimations.
\end{abstract}

 

\maketitle


\section{Introduction}

The proof of the uniqueness of weak solutions to cross-diffusion systems is a delicate and largely unsolved problem. When the diffusion equations exhibit specific structural properties, the uniqueness of weak solutions can be shown \cite{CLR17,GeJu18,ZaJu17}. The weak--strong uniqueness property represents a weaker uniqueness concept and was studied in, e.g., \cite{BBEP20,ChJu19-1,Hop22}. Another simplification is the analysis of cross-diffusion equations with a {\em triangular} diffusion matrix \cite{BDD25,ChDu25}. In this paper, we establish the uniqueness of bounded weak solutions for a class of two-species triangular cross-diffusion systems, which are motivated by nutrient-taxis models \cite{LMP13}. 

To encompass a broad range of models studied in the literature, we consider a general class of triangular two-species systems:
\begin{align}
  \pa_t u &= \diver\big(A_{11}(u,v)\na u + A_{12}(u,v)\na v\big)
  + R_1(u,v), \label{1.u} \\
  \pa_t v &= \diver\big(A_{22}(u,v)\na v\big) + R_2(u,v)
  \quad\mbox{in }\Omega,\ t>0, \label{1.v}
\end{align}
together with initial and no-flux boundary conditions:
\begin{align}
  & u(0) = u^0, \quad v(0) = v^0 \quad\mbox{in }\Omega, \label{1.ic} \\
  & \big(A_{11}(u,v)\na u + A_{12}(u,v)\na v\big)\cdot\nu = 0,\
  A_{22}(u,v)\na v\cdot\nu = 0 \quad\mbox{on }\pa\Omega,
  \label{1.bc}
\end{align}
where $\Omega\subset\R^d$ ($d\ge 1$) and $\nu$ is the exterior unit normal vector to $\pa\Omega$ (which is assumed to exist). The functions $A_{ij}(u,v)$ are the diffusion coefficients and $R_i(u,v)$ are reaction rates. A typical example is the motion of food-consuming bacteria in a nutrient-poor environment:
\begin{align}\label{1.taxis}
  \pa_t u = \diver(u^\alpha v\na u - \chi u^\beta v\na v) +f(u,v), \quad 
  \pa_t v = \Delta v - uv,
\end{align}
where $u$ and $v$ represent the density of the bacterial cells and the concentration of the nutrients, respectively, $\alpha\geq0$, $\beta>0$, $\chi>0$ are parameters, and $f(u,v)$ describes the effect of consumption. The model with $\alpha=1$, $\beta=2$ and $f(u,v)=uv$ was proposed in \cite{LMP13} and derived from a kinetic transport equation in the diffusion limit in \cite{Pla19}. The existence of global weak solutions with these parameters was proved in \cite{Win21,ZhLi2405} in one and two space dimensions. The boundedness of solutions was shown in \cite{ZhLi2405} in two dimensions, i.e., the constructed solution satisfies $\|u(t)\|_{L^\infty(\Omega)} + \|v(t)\|_{W^{1,\infty}(\Omega)} \le C$ for all $t>0$. In the one-dimensional setting, such solutions turn out to be unique when $1/u^0\in L^1(\Omega)$ \cite{ChDu25}. Surprisingly, the question of the uniqueness of weak solutions in the multi-dimensional setting or for other parameters $(\alpha,\beta)$ is open. In this paper, we fill this gap. 

For our result, we need some structural conditions on the diffusion coefficients $A_{ij}$ that include equations of the type \eqref{1.taxis}. We impose a condition that is slightly stronger than Lipschitz continuity. More precisely, for $\gamma\ge 1$, we call a function $f:[0,\infty)^2\to\R$ {\em locally $\gamma$-Lipschitz continuous} if for any $\rho>0$, there exists $C(\rho)>0$ such that for any $(u_1,v_1)$, $(u_2,v_2)\in[0,\rho]^2$,
\begin{align*}
  |f(u_1,v_1)-f(u_2,v_2)|\le C(\rho)
  \big(|u_1^\gamma-u_2^\gamma| + |v_1-v_2|\big).
\end{align*}
For instance, $f(u,v)=u^\gamma v$ with $\gamma\ge 1$ is locally $\gamma$-Lipschitz continuous. Any locally $\gamma$-Lipschitz continuous function is locally $\beta$-Lipschitz continuous for any $1\le\beta<\gamma$. 

We impose the following assumptions:
\begin{itemize}
\item[(A1)] Domain: $\Omega\subset\R^d$ ($d\ge 1$) is a bounded domain with smooth boundary.
\item[(A2)] Initial data: $u^0$, $v^0\in W^{1,\infty}(\Omega)$ such that $u^0\not\equiv 0$ is nonnegative and $v^0$ is positive in $\overline\Omega$.
\item[(A3)] Diffusion coefficients I: $A_{11}(u,v) = p(v)u^\alpha$ with $p:[0,\infty)\to\R$, $p\in C^2((0,\infty))$, $p>0$ in $(0,\infty)$, and $\alpha\ge 0$; $A_{22}(u,v)\ge a(v)$ with $a:[0,\infty)\to\R$, $a\in C^0((0,\infty))$, $a>0$ in $(0,\infty)$. 
\item[(A4)] Diffusion coefficients II: $A_{12}$, $A_{22}$ are locally $(\alpha/2+1)$-Lipschitz continuous.
\item[(A5)] Reactions: $R_i(u,v)=uq_i(v) + R_i^0(u,v)$ with $q_i:[0,\infty)\to\R$, $q_i\in C^1((0,\infty))$ and $R_i^0$ is locally $(\alpha/2+1)$-Lipschitz continuous ($i=1,2$).
\end{itemize}

We observe that the coefficients $A_{11}$ and $A_{22}$ may vanish if $v=0$ and $A_{11}$ vanishes if $u=0$ (in case $\alpha>0$). Thus, equation \eqref{1.u} is doubly degenerated with respect to $A_{11}$. A non-standard degeneracy was treated for some cross-diffusion systems in \cite{GeJu18} by exploiting additional a priori bounds, which are not available here. In the nutrient-taxis model \eqref{1.taxis}, we can argue in a different way: If an $L^\infty(\Omega)$ bound for $u$ is available and $v^0$ is strictly positive, the minimum principle for \eqref{1.v} provides a positive lower bound $v(t)\ge c>0$ in $\Omega$ \cite[Lemma 3.8]{Win24}. In this situation, the degeneracy of $A_{11}$ at $v=0$ is ruled out, but equation \eqref{1.u} is still degenerate with respect to $u$.

An example satisfying Assumptions (A3)--(A5) is given by equations \eqref{1.taxis} with $f(u,v)=uv$, i.e.
\begin{align*}
  A_{11} = u^\alpha v, \quad A_{12} = -\chi u^\beta v, \quad
  A_{22} = 1, \quad q_1 = v, \quad q_2 = -v, \quad R_1^0 = R_2^0 = 0,
\end{align*}
where $\beta\ge\alpha/2+1$. In particular, the choice $\alpha=1$ and $\beta=2$ is admissible.

Our definition of bounded weak solution is as follows.

\begin{definition}[Bounded weak solution]\label{def.weak}
We call $(u,v)$ a {\em global bounded weak solution} to \eqref{1.u}--\eqref{1.bc} if $u(t)\ge 0$ in $\Omega$ for $t>0$, $v(t)\ge c(T)>0$ in $\Omega$ for $0<t<T$ and any $T>0$,
\begin{align*}
  u\in L^\infty_{\rm loc}(0,\infty;L^\infty(\Omega))\cap
  L^2_{\rm loc}(0,\infty;H^1(\Omega)), \quad
  v\in L^\infty_{\rm loc}(0,\infty;W^{1,\infty}(\Omega)),
\end{align*}
it holds for all $\phi\in C_0^\infty(\overline\Omega\times[0,\infty))$ that
\begin{align}
  -\int_0^\infty&\int_\Omega u\pa_t\phi dxdt 
  - \int_\Omega u^0\phi(0)dx \label{1.weaku} \\
  &= -\int_0^\infty\int_\Omega
  \big(A_{11}(u,v)\na u+A_{12}(u,v)\na v\big)
  \cdot\na\phi dxdt
  + \int_0^\infty\int_\Omega R_1(u,v)\phi dxdt, \nonumber \\
  -\int_0^\infty&\int_\Omega v\pa_t\phi dxdt 
  - \int_\Omega v^0\phi(0)dx \label{1.weakv} \\
  &= -\int_0^\infty\int_\Omega A_{22}(u,v)\na v\cdot\na\phi dxdt
  + \int_0^\infty\int_\Omega R_2(u,v)\phi dxdt, \nonumber 
\end{align}
and the initial data \eqref{1.ic} is satisfied in the sense of $H^1(\Omega)'$.
\end{definition}

For degenerate problems, we cannot generally expect an $L^2(\Omega)$ gradient bound as required in Definition \ref{def.weak}. We claim that this condition is reasonable if $\alpha\in(0,2]$. Indeed, consider \eqref{1.taxis} with $\beta=2$ for illustration. Assume that $u\in L^\infty(\Omega_T)$, where $\Omega_T=\Omega\times(0,T)$. Then the minimum principle implies that $v\ge c>0$ in $\Omega_T$ for some $c>0$, and a formal computation shows that
\begin{align}\label{1.ei}
  \frac{d}{dt}\int_\Omega&\bigg(-\log u + \frac{\chi}{2}|\na v|^2\bigg)dx
  + \frac{4}{\alpha^2}\int_\Omega v|\na u^{\alpha/2}|^2 dx\\\nonumber
  &+ \chi\int_\Omega(|\Delta v|^2 + u|\na v|^2)dx
  + \int_\Omega vdx = 0.
\end{align}
Because of $v\ge c>0$, this yields a bound for $\na u^{\alpha/2}$ in $L^2(\Omega_T)$ and consequently $\na u = (2/\alpha)u^{1-\alpha/2}\na u^{\alpha/2}\in L^2(\Omega_T)$ if $\alpha\le 2$. Thus, the condition $u\in L^2(0,T;H^1(\Omega))$ in Definition \ref{def.weak} can be realized in this (and other) situations.

Our main result is as follows.

\begin{theorem}[Continuous dependence on the initial data \& uniqueness]
\label{thm.cont}
Let $T>0$ and let $(u_i,v_i)$ be a global bounded weak solution to \eqref{1.u}--\eqref{1.bc} in the sense of Definition \ref{def.weak} with the initial datum $(u_i^0,v_i^0)$ for $i=1,2$. Then there exists $C(T)>0$ depending on $T$ such that for $0<t<T$,
\begin{align}\label{1.cont}
  \|(u_1-u_2)(t)\|_{H^1(\Omega)'}^2 &+ \|(v_1-v_2)(t)\|_{L^2(\Omega)}^2
  + \bigg(\int_\Omega(u_1-u_2)(t)dx\bigg)^2 \\
  &\le C(T)\big(\|u_1^0-u_2^0\|_{L^2(\Omega)}^2 
  + \|v_1^0-v_2^0\|_{L^2(\Omega)}^2\big). \nonumber 
\end{align}
In particular, if $(u_1^0,v_1^0)=(u_2^0,v_2^0)$ then $(u_1,v_1)(t)=(u_2,v_2)(t)$ in $\Omega$, $t>0$, i.e., the solution is unique.
\end{theorem}

The proof is based on the $H^{-1}(\Omega)$ method. More precisely, let $(u_1,v_1)$ and $(u_2,v_2)$ be two bounded weak solutions to \eqref{1.u}--\eqref{1.bc} and let $\psi_i$ be the unique solution to $-\Delta\psi_i = u_i - \operatorname{meas}(\Omega)^{-1}\int_\Omega udx$ with homogeneous Neumann boundary conditions and $\int_\Omega\psi_i dx=0$ for $i=1,2$. The proof of Theorem \ref{thm.cont} is carried out in four steps:

\begin{enumerate}[1.]
\item Choose the test function $\phi\equiv 1$ in the weak formulation of the difference of equation \eqref{1.u} for $u_1$ and $u_2$. We obtain after some estimations, detailed in Lemma \ref{lem.ineq1},
\begin{align*}
  \bigg(\int_\Omega&(u_1-u_2)(t)dx\bigg)^2
  \le 4\bigg(\int_\Omega(u_1^0-u_2^0)dx\bigg)^2\\
  &+C_1\int_0^t\int_\Omega(u_1^{\alpha+1}-u_2^{\alpha+1})(u_1-u_2)dxds+ C_1 \int_0^tGds,
\end{align*}
where $C_1>0$ depends on $T$ but not on the solution and
\begin{align*}
  G := \bigg(\int_\Omega(u_1-u_2)dx\bigg)^2
  + \|\na(\psi_1-\psi_2)\|_{L^2(\Omega)}^2
  + \|v_1-v_2\|_{L^2(\Omega)}^2.
\end{align*}
\item Choose the test function $\psi_1-\psi_2$ in weak formulation of the difference of equation \eqref{1.u} for $u_1$ and $u_2$. In some sense, this test function introduces an inverse Laplacian which cancels with the ``Laplacian part'' of the term $\diver(A_{11}\na u) = (\alpha+1)^{-1}p(v)\Delta u^{\alpha+1} + \na A_{11}\cdot\na u$, yielding the second term on the left-hand side of the estimate
\begin{align*}
  \frac{d}{dt}\|\na(\psi_1-\psi_2)\|_{L^2(\Omega)}^2
  &+ c_0\int_\Omega(u_1^{\alpha+1}-u_2^{\alpha+1})(u_1-u_2)dx \\
  &\le \delta\|\na(v_1-v_2)\|_{L^2(\Omega)}^2 + C_2(\delta) G,
\end{align*}
where $\delta\in (0,1)$ is arbitrary and $c_0>0$, $C_2(\delta)>0$ are some constants (see Lemma \ref{lem.ineq2}). The second term on the left-hand side allows us to cancel the first term on the right-hand side of the inequality of Step 1.
\item Choose the test function $v_1-v_2$ in the weak formulation of the difference of equation \eqref{1.v} for $v_1$ and $v_2$ and estimate as follows (see lemma \ref{lem.ineq3}):
\begin{align*}
  \frac{d}{dt}\|v_1-v_2\|_{L^2(\Omega)}^2
  + c_1\|\na(v_1-v_2)\|_{L^2(\Omega)}^2
  \le C_3\int_\Omega(u_1^{\alpha+1}-u_2^{\alpha+1})(u_1-u_2)dx + C_3 G,
\end{align*}
where $C_3>0$ is another constant. The gradient term $\na(v_1-v_2)$ is a consequence of the condition $A_{22}>0$ in Assumption (A4). It allows us to eliminate the first term on the right-hand side of the inequality of Step 2.
\item Combine the three inequalities and choose $\delta$ in such a way that the terms involving $(u_1^{\alpha+1}-u_2^{\alpha+1})(u_1-u_2)$ and $\na(v_1-v_2)$ on the right-hand sides of the inequalities are eliminated. This gives $G(t)\leq C\int_0^tG(s)ds+CG(0)$ for some $C>0$, and an application of Gronwall's lemma finishes the proof.
\end{enumerate}

The uniqueness result can be applied to the double degenerate nutrient-taxis system \eqref{1.taxis}.

\begin{corollary}\label{coro}
Let $f:[0,\infty)^2\to\R$ be locally $(\alpha/2+1)$-Lipschitz continuous and let $\beta\ge \alpha/2+1$. Then any global bounded weak solution to the nutrient-taxis model \eqref{1.ic}--\eqref{1.taxis} is unique.
\end{corollary}

Table \ref{table} illustrates some variants of the nutrient-taxis system satisfying the assumptions of Corollary \ref{coro} and for which the existence of global bounded weak solutions was proved in the literature. 

\renewcommand{\arraystretch}{1.1}
\begin{table}[ht]
\begin{tabular}{|c|c|c|c|c|}
\hline
$\alpha$ & $\beta$ & $d$ & $f(u,v)$ & reference \\ \hline
1 & 2 & 1 & $uv$ & \cite{Win21} \\
1 & 2 & 2 & $k uv$ & \cite{ZhLi2405} \\
1 & $[3/2,2)$ & 2 & $kuv$ & \cite{Win24} \\ 
1 & $[3/2,2)$ & 3 & $kuv$ & \cite{TrYa25} \\ 
1 & 2 & 2 & $u(1-u)$ & \cite{LiWi25} \\
1 & 2 & $\ge 2$ & $u(\rho-\mu u^\kappa)$
& \cite{Pan24} \\
$[0,2]$ & $\alpha+1$ & 2 & $u(1-u)$ & \cite{ZhLi2410} \\
$[0,2]$ & $\alpha+1$ & 2 & $uv$ & \cite{ZhLi2411} \\
\hline
\end{tabular}
\medskip
\caption{Parameter values for which the existence of global bounded weak solutions to \eqref{1.ic}--\eqref{1.bc} is proved. In the table, the parameters are $k\ge 0$, $\rho>0$, $\mu>0$, $\kappa>d/2$.}
\label{table}
\end{table}

In all cases of Table \ref{table}, the regularity $u\in L^\infty_{\rm loc}(0,\infty;L^\infty(\Omega))$, $v\in L^\infty_{\rm loc}(0,\infty;W^{1,\infty}(\Omega))$ has been shown. We only need to verify that $\na u\in L^2(\Omega_T)$ for all $T>0$. In fact, it has been shown in the references in Table \ref{table} that $\|u\|_{L^\infty(\Omega_T)}\le K$ for some $K>0$. Then, by the minimum principle, $v\ge k>0$ in $\Omega_T$ for some $k>0$. In the literature, entropy inequalities related to the entropy $\int_\Omega(h(u)+\frac12|\na v|^2)dx$ with $h(u)=-\log u$ or $h(u)=u(\log u-1)$ have been derived, and the entropy production includes the integrals
\begin{equation}\label{1.ei2} 
\begin{aligned}
  \int_0^T\int_\Omega |\na u|^2 dxdt
  &\le \frac{K}{k}\int_0^T\int_\Omega\frac{v}{u}|\na u|^2 dxdt \le C
  \quad\mbox{or} \\
  \int_0^T\int_\Omega |\na u|^2 dxdt
  &\le \frac{1}{k}\int_0^T\int_\Omega v|\na u|^2 dxdt;
\end{aligned}
\end{equation}
also see \eqref{1.ei}. This yields the desired $L^2(\Omega_T)$ bound for $\na u$. 

\begin{remark}[Existence of solutions]\rm
The existence of global bounded weak solutions to the special class of systems \eqref{1.taxis} was shown in the literature in some special settings; see Table \ref{table}. An existence result for the general model \eqref{1.u}--\eqref{1.bc} is open. One delicate point is the proof of the boundedness of $u$. If this property can be proved, the maximum principle implies that $v$ is bounded in any finite time interval, and parabolic regularity theory leads to $v\in L^r(0,T;W^{2,r}(\Omega))\hookrightarrow L^r(0,T;W^{1,\infty}(\Omega))$ for $r>d$. If the corresponding estimate is independent of $r$, the limit $r\to\infty$ gives $v\in L^\infty(0,T;W^{1,\infty}(\Omega))$. The gradient bound for $u$ follows if suitable entropy estimates are available (see, e.g., \eqref{1.ei2}). Such an estimate requires a particular structure of the reaction rate $R_2$. We conclude that the existence of solutions in the sense of Definition \ref{def.weak} strongly relies on the structure of the nonlinearities $A_{ij}$ and $R_i$ and a general result seems to be out of reach.
\qed\end{remark}

The paper is organized as follows. Some auxiliary results are presented in Section \ref{sec.aux}, and Section \ref{sec.proof} is devoted to the proof of Theorem \ref{thm.cont}. 


\section{Auxiliary results}\label{sec.aux}

Let $T>0$ and $u_i\in L^2(\Omega_T)$ for $i=1,2$, recalling that $\Omega_T:=\Omega\times(0,T)$. There exists a unique solution $\psi_i\in L^2(0,T;H^2(\Omega))$ to
\begin{align}\label{2.psi}
  -\Delta\psi_i = u_i - \fint_\Omega u_i dx\quad\mbox{in }\Omega,
  \quad \na\psi_i\cdot\nu = 0\quad\mbox{on }\pa\Omega, \quad
  \int_\Omega \psi_i dx = 0,
\end{align}
where $\fint u_i dx= \operatorname{meas}(\Omega)^{-1}\int_\Omega u_i dx$. If additionally $\pa_t u_i\in L^2(0,T;H^1(\Omega)')$, we have $\pa_t\psi_i\in L^2(0,T;H^1(\Omega))$. In particular, $t\mapsto\frac12\|\na(\psi_1-\psi_2)\|_{L^2(\Omega)}^2$ is absolutely continuous, which yields \cite[Sec.~5.9.2, Theorem 3]{Eva02}
\begin{align}\label{2.psi12}
  \frac12\|&\na(\psi_1-\psi_2)(t)\|_{L^2(\Omega)}^2
  - \frac12\|\na(\psi_1-\psi_2)(0)\|_{L^2(\Omega)}^2 \\
  &= \frac12\int_0^t\frac{d}{dt}\|\na(\psi_1-\psi_2)
  \|_{L^2(\Omega)}^2 ds
  = \int_0^t\langle\pa_t\na(\psi_1-\psi_2),
  \na(\psi_1-\psi_2)\rangle ds \nonumber \\
  &= \int_0^t\langle\pa_t(u_1-u_2),\psi_1-\psi_2\rangle ds, \nonumber 
\end{align}
where $\langle\cdot,\cdot\rangle$ is the dual product between $H^1(\Omega)'$ and $H^1(\Omega)$. The weak formulations \eqref{1.weaku} and \eqref{1.weakv} holds for test functions $\phi\in C_0^\infty(\overline\Omega\times[0,\infty))$. We wish to allow for test functions in $C^\infty(\overline\Omega\times[0,t])$ for $t\in(0,T)$. The following result was proved in step 1 of the proof of \cite[Lemma 13]{ChJu19}.

\begin{lemma}\label{lem.weak1}
Let $w$ be a weak solution to $\pa_t w + \diver F=g$ in $\Omega$, $t>0$, for integrable functions $F$ and $g$ in the sense
\begin{align}\label{2.weak}
  -\int_0^T\int_\Omega w\pa_t\chi dxdt - \int_\Omega w(0)\chi(0)dx
  = \int_0^T\int_\Omega(F\cdot\na\chi + g\chi)dxdt
\end{align}
for all $\chi\in C_0^\infty(\overline\Omega\times[0,T))$. Then, for any $t\in[0,T]$ and $\phi\in C^\infty(\overline\Omega\times[0,T])$,
\begin{align*}
  -\int_0^t\int_\Omega w\pa_s\phi dxds 
  + \int_\Omega\big(w(t)\phi(t)-w(0)\phi(0)\big)dx
  = \int_0^t\int_\Omega(F\cdot\na\chi + g\chi)dxds.
\end{align*}
\end{lemma}

If $F$ and $g$ are square integrable, the weak formulation can be slightly strengthened.

\begin{lemma}\label{lem.weak2}
Let $w$ be a weak solution to $\pa_t w + \diver F=g$ in $\Omega$, $t>0$, for $F$, $g\in L^2(\Omega_T)$ in the sense of \eqref{2.weak} for all $\chi\in C_0^\infty(\overline\Omega\times[0,T))$. Then $\pa_t w\in L^2(0,T;H^1(\Omega)')$ and for $t\in[0,T]$ and $\phi\in L^2(0,T;H^1(\Omega))$,
\begin{align*}
  \int_0^t\langle\pa_s w,\phi\rangle ds 
  = \int_0^t\int_\Omega(F\cdot\na\chi + g\chi)dxds.
\end{align*}
\end{lemma}

We use some elementary inequalities.

\begin{lemma}
Let $\alpha\ge 0$. There exists $C(\alpha)>0$ such that for all $u_1$, $u_2\ge 0$,
\begin{align}\label{alpha1}
  \big(u_1^{\alpha/2+1}-u_2^{\alpha/2+1}\big)^2
  &\le C(\alpha)(u_1^{\alpha+1}-u_2^{\alpha+1})(u_1-u_2).
\end{align}
Let $M>0$. Then there exists $C(M)>0$ such that for all $0\le u_1\le M$, $0\le u_2\le M$,
\begin{align}\label{alpha2}
  (u_1^{\alpha+1}-u_2^{\alpha+1})^2 
  \le C(M)(u_1^{\alpha+1}-u_2^{\alpha+1})(u_1-u_2).
\end{align}
\end{lemma}

\begin{proof}
Inequality \eqref{alpha1} is proved in \cite[Lemma A.3]{CJS16} with $C(\alpha) = (\alpha+2)^2/(4\alpha+4)$. By the mean-value theorem, for some $\xi$ between $u_1$ and $u_2$,
\begin{align*}
  (u_1^{\alpha+1}-u_2^{\alpha+1})^2
  &= (\alpha+1)\xi^\alpha(u_1-u_2)(u_1^{\alpha+1}-u_2^{\alpha+1}) \\
  &\le (\alpha+1)M^\alpha(u_1^{\alpha+1}-u_2^{\alpha+1})(u_1-u_2).
\end{align*}
This shows inequality \eqref{alpha2} with $C(M)=(\alpha+1)M^\alpha$.
\end{proof}


\section{Proof of Theorem \ref{thm.cont}}\label{sec.proof}

Let $T>0$ and let $(u_i,v_i)$ be a global bounded weak solution to \eqref{1.u}--\eqref{1.bc} in the sense of Definition \ref{def.weak}. The proof of Theorem \ref{thm.cont} is based on several estimations relating some norms of the differences $u_1-u_2$, $\psi_1-\psi_2$, and $v_1-v_2$.

\begin{lemma}\label{lem.ineq1}
There exists a constant $C_1>0$ depending on $T$ such that for $0<t<T$,
\begin{align*}
  \bigg(\int_\Omega(u_1-u_2)(t)dx\bigg)^2 
  &\le 4\bigg(\int_\Omega(u_1^0-u_2^0)dx\bigg)^2 
  + C_1\int_0^t\bigg(\int_\Omega(u_1-u_2)dx\bigg)^2ds \\
  &\phantom{xx}+ C_1\int_0^t\|\na(\psi_1-\psi_2)\|_{L^2(\Omega)}^2ds
  + C_1\int_0^t\|v_1-v_2\|_{L^2(\Omega)}^2ds \\
  &\phantom{xx}+ C_1\int_0^t\int_\Omega(u_1^{\alpha+1}-u_2^{\alpha+1})
  (u_1-u_2)dxds.
\end{align*}
\end{lemma}

\begin{proof}
According to Lemma \ref{lem.weak1}, the weak formulation \eqref{1.weaku} can be written as
\begin{align*}
  -&\int_0^t\int_\Omega u_i\pa_s\phi dxds
  + \int_\Omega\big(u_i(t)\phi(t)-u_i(0)\phi(0)\big)dx \\
  &= -\int_0^t\int_\Omega\big(A_{11}(u_i,v_i)\na u_i
  + A_{12}(u_i,v_i)\na v_i\big)\cdot\na\phi dxds
  + \int_0^t\int_\Omega R_1(u_i,v_i)\phi dxds
\end{align*}
for $i=1,2$, $t\in[0,T]$, and $\phi\in C^\infty(\overline\Omega\times[0,T])$. This becomes for $\phi\equiv 1$:
\begin{align*}
  \int_\Omega(u_i(t)-u_i(0))dx
  = \int_0^t\int_\Omega R_1(u_i,v_i)\phi dxds, \quad i=1,2.
\end{align*}
We take the difference of this equation for $i=1$ and $i=2$, insert the definition of $R_i$ according to Assumption (A5), and take the square:
\begin{align}\label{3.ineq1}
  & \bigg(\int_\Omega(u_1-u_2)(t)dx\bigg)^2 
  \le 4\bigg(\int_\Omega(u_1^0-u_2^0)dx\bigg)^2
  + I_1 + I_2 + I_3, \quad\mbox{where} \\
  & I_1 = 4\bigg(\int_0^t\int_\Omega(u_1-u_2)q_1(v_1)dxds\bigg)^2, 
  \nonumber \\
  & I_2 = 4\bigg(\int_0^t\int_\Omega u_2(q_1(v_1)-q_1(v_2))
  dxds\bigg)^2, \nonumber \\
  & I_3 = 4\bigg(\int_0^t\int_\Omega(R_1^0(u_1,v_1)-R_2^0(u_2,v_2))
  dxds\bigg)^2. \nonumber 
\end{align}
To estimate $I_1$, we apply the Cauchy--Schwarz inequality, insert \eqref{2.psi} for $u_1-u_2$, and integrate by parts:
\begin{align*}
  I_1 &\le 4t\int_0^t\bigg(\int_\Omega(u_1-u_2)q_1(v_1)dx\bigg)^2 ds \\
  &= 4t\int_0^t\bigg(\int_\Omega q_1'(v_1)\na(\psi_1-\psi_2)
  \cdot\na v_1 dx + \fint_\Omega(u_1-u_2)dx
  \int_\Omega q_1(v_1)dx\bigg)^2 ds \\
  &\le C\int_0^t\int_\Omega|\na(\psi_1-\psi_2)|^2 dxds
  + C\int_0^t\bigg(\int_\Omega(u_1-u_2)dx\bigg)^2 ds,
\end{align*} 
where $C>0$ depends on $T$ and on the $L^\infty(0,T;W^{1,\infty}(\Omega))$ norm of $v_1$. Next, we use the boundedness of $v_2$ to find that
\begin{align*}
  I_2 &\le C\int_0^t\int_\Omega(q_1(v_1)-q_1(v_2))^2 dxds \\
  &= C\int_0^t\int_\Omega\bigg(\int_0^1
  q_1'(\theta v_1+(1-\theta)v_2)d\theta\bigg)^2(v_1-v_2)^2 dxds \\
  &\le C\int_0^t\int_\Omega(v_1-v_2)^2 dxds,
\end{align*}
where here and in the following, the value of the constant $C>0$ changes from line to line. It follows from the $(\alpha/2+1)$-Lipschitz continuity of $R_1^0$ that
\begin{align*}
  I_3 &\le C\int_0^t\int_\Omega
  \big(R_1^0(u_1,v_1)-R_1^0(u_2,v_2)\big)^2 dxds \\
  &\le C\int_0^t\int_\Omega\big((u_1^{\alpha/2+1}-u_2^{\alpha/2+1})^2
  + (v_1-v_2)^2\big)dxds. 
\end{align*}
To estimate the term involving $u_i^{\alpha/2+1}$, we use inequality \eqref{alpha1}, which yields
\begin{align}\label{3.I3}
  I_3 \le C\int_0^t\int_\Omega
  (u_1^{\alpha+1}-u_2^{\alpha+1})(u_1-u_2)dxds
  + C\int_0^t\int_\Omega(v_1-v_2)^2 dxds.
\end{align}
Inserting the estimates for $I_1$, $I_2$, and $I_3$ into \eqref{3.ineq1} finishes the proof.
\end{proof}

The integral over $(u_1^{\alpha+1}-u_2^{\alpha+1})(u_1-u_2)$ can be controlled up to a small contribution of the $L^2(\Omega)$ norm of $\na(v_1-v_2)$, as shown in the following lemma.

\begin{lemma}\label{lem.ineq2}
Let $\delta\in (0,1)$. Then there exist constants $c_0>0$ and $C_2>0$ depending on $T$ such that
\begin{align*}
  \|\na&(\psi_1-\psi_2)(t)\|_{L^2(\Omega)}^2
  + c_0\int_0^t\int_\Omega(u_1^{\alpha+1}-u_2^{\alpha+1})
  (u_1-u_2)dxds \\
  &\le \|\na(\psi_1-\psi_2)(0)\|_{L^2(\Omega)}^2
  + C_2\int_0^t\bigg(\int_\Omega(u_1-u_2)dx\bigg)^2 ds \\
  &\phantom{xx}
  + C_2\delta^{-1}\int_0^t\|\na(\psi_1-\psi_2)\|_{L^2(\Omega)}^2 ds \\
  &\phantom{xx}+ C_2\int_0^t\|v_1-v_2\|_{L^2(\Omega)}^2 ds
  + \delta\int_0^t\|\na(v_1-v_2)\|_{L^2(\Omega)}^2 ds.
\end{align*}
\end{lemma}

\begin{proof}
We split the proof into several steps.

{\em Step 1: Preparations.} According to Lemma \ref{lem.weak2}, the weak formulation of the difference $u_1-u_2$ can be written for all $\phi\in L^2(0,T;H^1(\Omega))$ as
\begin{align*}
  \int_0^t\langle\pa_s(u_1-u_2),\phi\rangle ds
  &= -\int_0^t\int_\Omega\big(A_{11}(u_1,v_1)\na u_1
  - A_{11}(u_2,v_2)\na u_2\big)\cdot\na\phi dxds \\
  &\phantom{xx}- \int_0^t\int_\Omega\big(A_{12}(u_1,v_1)\na v_1
  - A_{12}(u_2,v_2)\na v_2\big)\cdot\na\phi dxds \\
  &\phantom{xx}+ \int_0^t\int_\Omega
  \big(R_1(u_1,v_1)-R_1(u_2,v_2)\big)\phi dxds.
\end{align*}
This gives, with the test function $\phi=\psi_1-\psi_2$ and identity \eqref{2.psi12},
\begin{align*}
  & \frac12\int_\Omega|\na(\psi_1-\psi_2)(t)|^2 dx
  - \frac12\int_\Omega|\na(\psi_1-\psi_2)(0)|^2 dx
  = I_4+I_5+I_6, \quad\mbox{where} \\
  & I_4 = -\int_0^t\int_\Omega\big(A_{11}(u_1,v_1)\na u_1
  - A_{11}(u_2,v_2)\na u_2\big)\cdot\na(\psi_1-\psi_2) dxds, \\
  & I_5 = -\int_0^t\int_\Omega\big(A_{12}(u_1,v_1)\na v_1
  - A_{12}(u_2,v_2)\na v_2\big)\cdot\na(\psi_1-\psi_2) dxds, \\
  & I_6 = \int_0^t\int_\Omega\big(R_1(u_1,v_1)-R_1(u_2,v_2)\big)
  (\psi_1-\psi_2)dxds.
\end{align*}

{\em Step 2: Estimation of $I_4$.} We rewrite the term $I_4$ by inserting Assumption (A3):
\begin{align*}
  I_4 
  &= -\frac{1}{\alpha+1}\int_0^t\int_\Omega
  p(v_1)\na(u_1^{\alpha+1}-u_2^{\alpha+1})\cdot\na(\psi_1-\psi_2)dxds \\
  &\phantom{xx} - \frac{1}{\alpha+1}\int_0^t\int_\Omega  
  (p(v_1)-p(v_2))\na u_2^{\alpha+1}\cdot\na(\psi_1-\psi_2) dxds.
\end{align*}
To get rid of the gradient in $\na(u_1^{\alpha+1}-u_2^{\alpha+1})$, we integrate by parts. Then $I_4=I_{41}+\cdots+I_{44}$, where
\begin{align*}
  I_{41} &= \frac{1}{\alpha+1}\int_0^t\int_\Omega
  (u_1^{\alpha+1}-u_2^{\alpha+1})\na p(v_1)\cdot\na(\psi_1-\psi_2) 
  dxds, \\
  I_{42} &= \frac{1}{\alpha+1}\int_0^t\int_\Omega p(v_1)
  (u_1^{\alpha+1}-u_2^{\alpha+1})\Delta(\psi_1-\psi_2) dxds, \\
  I_{43} &= \frac{1}{\alpha+1}\int_0^t\int_\Omega
  u_2^{\alpha+1}\na(p(v_1)-p(v_2))\cdot\na(\psi_1-\psi_2) dxds, \\
  I_{44} &=  \frac{1}{\alpha+1}\int_0^t\int_\Omega(p(v_1)-p(v_2))
  u_2^{\alpha+1}\Delta(\psi_1-\psi_2)dxds.
\end{align*}

We start with $I_{42}$. For this, we insert \eqref{2.psi} and use H\"older's inequality:
\begin{align*}
  I_{42} &= -\frac{1}{\alpha+1}\int_0^t\int_\Omega p(v_1)
  (u_1^{\alpha+1}-u_2^{\alpha+1})(u_1-u_2)dxds \\
  &\phantom{xx}+ \frac{1}{\alpha+1}\int_0^t\int_\Omega p(v_1)
  (u_1^{\alpha+1}-u_2^{\alpha+1})dx
  \bigg(\fint_\Omega(u_1-u_2)dx\bigg)ds \\
  &\le -c_0\int_0^t\int_\Omega
  (u_1^{\alpha+1}-u_2^{\alpha+1})(u_1-u_2)dxds \\
  &\phantom{xx}+ C\int_0^t\bigg(\int_\Omega
  (u_1^{\alpha+1}-u_2^{\alpha+1})^2dx\bigg)^{1/2}
  \bigg|\int_\Omega(u_1-u_2)dx\bigg|ds,
\end{align*}
where $c_0>0$ depends on the positive lower bound of $p(v_1)$ on $\Omega_T$. The last term is estimated by means of inequality \eqref{alpha2} and Young's inequality as
\begin{align*}
  C\int_0^t\bigg(&\int_\Omega
  (u_1^{\alpha+1}-u_2^{\alpha+1})^2dx\bigg)^{1/2}
  \bigg|\int_\Omega(u_1-u_2)dx\bigg|ds \\
  &\le C\int_0^t\bigg(\int_\Omega (u_1^{\alpha+1}-u_2^{\alpha+1})
  (u_1-u_2)dx\bigg)^{1/2}\bigg|\int_\Omega(u_1-u_2)dx\bigg|ds \\
  &\le \frac{c_0}{8}\int_0^t\int_\Omega 
  (u_1^{\alpha+1}-u_2^{\alpha+1})(u_1-u_2)dxds
  + C\int_0^t\bigg(\int_\Omega(u_1-u_2)dx\bigg)^2ds.
\end{align*}
This yields
\begin{align*}
  I_{42} \le -\frac{7c_0}{8}\int_0^t\int_\Omega
  (u_1^{\alpha+1}-u_2^{\alpha+1})(u_1-u_2)dxds
  + C\int_0^t\bigg(\int_\Omega(u_1-u_2)dx\bigg)^2ds.
\end{align*}

It follows from Young's inequality and \eqref{alpha2} that 
\begin{align*}
  I_{41} &\le C\int_0^t\int_\Omega|u_1^{\alpha+1}-u_2^{\alpha+1}|
  |\na(\psi_1-\psi_2)|dxds \\
  &\le \frac{c_0}{8}\int_0^t\int_\Omega
  (u_1^{\alpha+1}-u_2^{\alpha+1})(u_1-u_2)dxds
  + C\int_0^t\int_\Omega|\na(\psi_1-\psi_2)|^2dxds.
\end{align*}
For the estimate of $I_{43}$, we first deduce from
\begin{align*}
  p(v_1)-p(v_2) 
  = \int_0^1 p'(\theta v_1+(1-\theta)v_2)d\theta(v_1-v_2)
\end{align*}
that
\begin{align*}
  |\na(p(v_1)-p(v_2))| &= \bigg|\int_0^1 p''(\theta v_1+(1-\theta)v_2)
  (\theta\na v_1+(1-\theta)\na v_2)d\theta(v_1-v_2) \\
  &\phantom{xx}+ \int_0^1 p'(\theta v_1+(1-\theta)v_2)d\theta
  \na(v_1-v_2)\bigg| \\
  &\le C\big(|v_1-v_2| + |\na(v_1-v_2)|\big),
\end{align*}
where $C>0$ depends on the $L^\infty(\Omega_T)$ norm of $\na v_i$ for $i=1,2$. This shows that, for any $\delta\in(0,1)$,
\begin{align*}
  I_{43} &\le C\int_0^t\int_\Omega\big(|v_1-v_2| + |\na(v_1-v_2)|\big)
  |\na(\psi_1-\psi_2)|dxds \\
  &\le \frac{\delta}{4}\int_0^t\int_\Omega|\na(v_1-v_2)|^2 dxds
  + C\int_0^t\int_\Omega(v_1-v_2)^2dxds \\
  &\phantom{xx}+ C\delta^{-1}\int_0^t\int_\Omega
  |\na(\psi_1-\psi_2)|^2 dxds.
\end{align*}

Finally, we insert equation \eqref{2.psi} to find that
\begin{align}\label{3.I44}
  I_{44} &= -\frac{1}{\alpha+1}\int_0^t\int_\Omega
  (p(v_1)-p(v_2))u_2^{\alpha+1}(u_1-u_2)dxds \\
  &\phantom{xx}+ \frac{1}{\alpha+1}\int_0^t\int_\Omega
  (p(v_1)-p(v_2))u_2^{\alpha+1}dx
  \bigg(\fint_\Omega(u_1-u_2)dx\bigg)ds \nonumber \\
  &=: I_{441}+I_{442}. \nonumber 
\end{align}
Let $F(u_1,u_2):=(1+\alpha)\int_0^1(\theta u_1+(1-\theta) u_2)^\alpha
d\theta$. Then we infer from
\begin{align*}
  (u_1^{\alpha+1}-u_2^{\alpha+1})(u_1-u_2) = F(u_1,u_2)(u_1-u_2)^2
\end{align*}
and Young's inequality that
\begin{align*}
  I_{441} &\le \frac{1}{\alpha+1}\int_0^t\int_\Omega
  \mathrm{1}_{\{u_1+u_2>0\}}|p(v_1)-p(v_2)|
  \frac{u_2^{\alpha+1}}{\sqrt{F(u_1,u_2)}}
  \sqrt{F(u_1,u_2)(u_1-u_2)^2}dx ds \\
  &\le \frac{c_0}{8}\int_0^t\int_\Omega
  (u_1^{\alpha+1}-u_2^{\alpha+1})(u_1-u_2)dxds \\
  &\phantom{xx}+ C\int_0^t\int_\Omega
  \mathrm{1}_{\{u_1+u_2>0\}}(p(v_1)-p(v_2))^2
  \frac{u_2^{2\alpha+2}}{F(u_1,u_2)}dx ds.
\end{align*}
We claim that the quotient in the last integral is bounded. Indeed, this is true if $u_2=0$ (but $u_1>0$). Otherwise, if $u_2>0$, we have
\begin{align*}
  F(u_1,u_2) \ge (\alpha+1)\int_0^1(1-\theta)^\alpha u_2^{\alpha}
  d\theta = u_2^\alpha,
\end{align*}
which yields $u_2^{2\alpha+2}/F(u_1,u_2)\le u_2^{\alpha+2}\le C$. Therefore, by the Lipschitz continuity of $p$ on finite intervals,
\begin{align*}
  I_{441} \le \frac{c_0}{8}\int_0^t\int_\Omega
  (u_1^{\alpha+1}-u_2^{\alpha+1})(u_1-u_2)dxds
  + C\int_0^t\int_\Omega(v_1-v_2)^2 dxds.
\end{align*}
Applying H\"older's and Young's inequality to $I_{442}$, we obtain
\begin{align*}
  I_{442} &\le C\int_0^t\int_\Omega|v_1-v_2|dx
  \bigg|\int_\Omega(u_1-u_2)dx\bigg|ds \\
  &\le C\int_0^t\int_\Omega(v_1-v_2)^2 dxds 
  + C\int_0^t\bigg(\int_\Omega(u_1-u_2)dx\bigg)^2 ds.
\end{align*}
Summarizing the estimates for $I_{441}$ and $I_{442}$, we conclude from \eqref{3.I44} that
\begin{align*}
  I_{44} &\le \frac{c_0}{8}\int_0^t\int_\Omega
  (u_1^{\alpha+1}-u_2^{\alpha+1})(u_1-u_2)dxds
  + C\int_0^t\int_\Omega(v_1-v_2)^2 dxds \\
  &\phantom{xx}+ C\int_0^t\bigg(\int_\Omega(u_1-u_2)dx\bigg)^2 ds.
\end{align*}
The estimates for $I_{41},\ldots,I_{44}$ show that
\begin{align*}
  I_4 &\le -\frac{5c_0}{8}\int_0^t\int_\Omega
  (u_1^{\alpha+1}-u_2^{\alpha+1})(u_1-u_2)dxds
  + C\delta^{-1}\int_0^t\int_\Omega|\na(\psi_1-\psi_2)|^2 dxds \\
  &\phantom{xx}+ C\int_0^t\bigg(\int_\Omega(u_1-u_2)dx\bigg)^2ds
  + C\int_0^t\int_\Omega(v_1-v_2)^2 dxds \\
  &\phantom{xx}+ \frac{\delta}{4}\int_0^t\int_\Omega
  |\na(v_1-v_2)|^2 dxds.
\end{align*}

{\em Step 3: Estimation of $I_5$.} We write $I_5=I_{51}+I_{52}$, where
\begin{align*}
  I_{51} &= -\int_0^t \int_{\Omega}\big(A_{12}(u_1,v_1)
  - A_{12}(u_2,v_2)\big)\na v_1\cdot\na(\psi_1-\psi_2) dxds, \\
  I_{52} &= -\int_0^t \int_{\Omega} A_{12}(u_2,v_2)\na(v_1-v_2) 
  \cdot\na(\psi_1-\psi_2) dxds.
\end{align*}
We exploit the $(\alpha/2+1)$-Lipschitz continuity of $A_{12}$ and use inequality \eqref{alpha1}:
\begin{align*}
  I_{51} &\le C\int_0^t\int_\Omega
  \big(|u_1^{\alpha/2+1}-u_2^{\alpha/2+1}| + |v_1-v_2|\big)
  |\na(\psi_1-\psi_2)|dxds \\
  &\le C\int_0^t\int_\Omega\Big(\sqrt{(u_1^{\alpha+1}-u_2^{\alpha+1})
  (u_1-u_2)} + |v_1-v_2|\Big)|\na(\psi_1-\psi_2)|dxds \\
  &\le \frac{c_0}{8}\int_0^t\int_\Omega(u_1^{\alpha+1}-u_2^{\alpha+1})
  (u_1-u_2)dxds
  + C\int_0^t\int_\Omega(v_1-v_2)^2 dxds \\
  &\phantom{xx}+ C\int_0^t\int_\Omega|\na(\psi_1-\psi_2)|^2 dxds.
\end{align*}
Thus, in view of
\begin{align*}
  I_{52} \le \frac{\delta}{4}\int_0^t\int_\Omega|\na(v_1-v_2)|^2 dxds
  + C\delta^{-1}\int_0^t\int_\Omega|\na(\psi_1-\psi_2)|^2 dxds,
\end{align*}
we infer that
\begin{align*}
  I_5 &\le \frac{c_0}{8}\int_0^t\int_\Omega(u_1^{\alpha+1}-u_2^{\alpha+1})
  (u_1-u_2)dxds 
  + C\delta^{-1}\int_0^t\int_\Omega|\na(\psi_1-\psi_2)|^2 dxds \\
  &\phantom{xx}+ C\int_0^t\int_\Omega(v_1-v_2)^2 dxds
  +  \frac{\delta}{4}\int_0^t\int_\Omega|\na(v_1-v_2)|^2 dxds.
\end{align*}

{\em Step 4: Estimation of $I_6$.} We use Assumption (A5) to formulate $I_6=I_{61}+I_{62}+I_{63}$, where
\begin{align*}
  I_{61} &= \int_0^t\int_\Omega(u_1-u_2)q_1(v_1)(\psi_1-\psi_2)dxds, \\
  I_{62} &= \int_0^t\int_\Omega u_2(q_1(v_1)-q_1(v_2))
  (\psi_1-\psi_2)dxds, \\
  I_{63} &= \int_0^t\int_\Omega
  \big(R_1^0(u_1,v_1)-R_1^0(u_2,v_2)\big)(\psi_1-\psi_2)dxds.
\end{align*}
We insert \eqref{2.psi} and integrate by parts:
\begin{align*}
  I_{61} &= -\int_0^t\int_\Omega\Delta(\psi_1-\psi_2)
  q_1(v_1)(\psi_1-\psi_2)dxds \\
  &\phantom{xx}+ \int_0^t\int_\Omega q_1(v_1)(\psi_1-\psi_2)dx
  \bigg(\fint_\Omega(u_1-u_2)dx\bigg)ds \\
  &= \int_0^t\int_\Omega\na(\psi_1-\psi_2)\cdot
  \big(q_1'(v_1)\na v_1(\psi_1-\psi_2) 
  + q_1(v_1)\na(\psi_1-\psi_2)\big)dxds \\
  &\phantom{xx}+ \int_0^t\int_\Omega q_1(v_1)(\psi_1-\psi_2)dx
  \bigg(\fint_\Omega(u_1-u_2)dx\bigg)ds \\
  &\le C\int_0^t\int_\Omega|\na(\psi_1-\psi_2)|^2 dxds
  + C\int_0^t\int_\Omega(\psi_1-\psi_2)^2 dxds \\
  &\phantom{xx}+ C\int_0^t\bigg(\int_\Omega(u_1-u_2)dx\bigg)^2 ds \\
  &\le C\int_0^t\int_\Omega|\na(\psi_1-\psi_2)|^2 dxds
  + C\int_0^t\bigg(\int_\Omega(u_1-u_2)dx\bigg)^2 ds,
\end{align*}
where we have applied the Poincar\'e inequality in the last step (this is possible since $\int_\Omega(\psi_1-\psi_2)dx=0$). Similarly, by the Lipschitz continuity of $q$ and Poincar\'e's inequality again,
\begin{align*}
  I_{62} \le C\int_0^t\int_{\Omega}(v_1 - v_2)^2 dxds
  + C\int_0^t\int_{\Omega} |\na(\psi_1 - \psi_2)|^2  dxds.
\end{align*}
It follows from estimate \eqref{3.I3}, for any $\eta>0$, that
\begin{align*}
  I_{63} &\le \eta\int_0^t\int_{\Omega}
  \big(R_1^0(u_1,v_1) - R_1^0(u_2,v_2)\big)^2dxds
  + C(\eta)\int_0^t\int_{\Omega} (\psi_1 - \psi_2)^2 dxds \\
  &\le C\eta\int_0^t \int_\Omega \big(u_1^{\alpha+1}-u_2^{\alpha+1}\big)(u_1-u_2)dxds
  + C\eta\int_0^t\int_\Omega(v_1-v_2)^2 dx ds \\
  &\phantom{xx}+ C(\eta)\int_0^t\int_\Omega|\na(\psi_1-\psi_2)|^2dxds.
\end{align*}
We choose $\eta>0$ such that $C\eta\le c_0/8$. Then, summarizing the estimates for $I_{61}$, $I_{62}$, and $I_{63}$,
\begin{align*}
  I_6 &\le \frac{c_0}{8}\int_0^t \int_\Omega \big(u_1^{\alpha+1}-u_2^{\alpha+1}\big)(u_1-u_2)dxds
  + C\int_0^t\int_\Omega|\na(\psi_1-\psi_2)|^2dxds \\
  &\phantom{xx}+ C\int_0^t\int_\Omega(v_1-v_2)^2 dx ds
  + C\int_0^t\bigg(\int_\Omega(u_1-u_2)dx\bigg)^2 ds.
\end{align*}
Putting together the inequalities for $I_4$, $I_5$, and $I_6$ concludes the proof. 
\end{proof}

It remains to find an estimate for the $H^1(\Omega)$ norm of $v_1-v_2$.

\begin{lemma}\label{lem.ineq3}
There exist constants $c_1>0$ and $C_3>0$ depending on $T$ such that
\begin{align*}
  \|(v_1&-v_2)(t)\|_{L^2(\Omega)}^2 
  + c_1\int_0^t\|\na(v_1-v_2)\|_{L^2(\Omega)}^2 ds
  \le \|(v_1-v_2)(0)\|_{L^2(\Omega)}^2 \\
  &+ C_3\int_0^t\int_\Omega\big(u_1^{1+\alpha}-u_2^{1+\alpha}\big)
  (u_1-u_2) dx ds
  + C_3\int_0^t\|\na(\psi_1-\psi_2)\|_{L^2(\Omega)}^2ds \\
  &+ C_3\int_0^t\bigg(\int_\Omega(u_1-u_2)dx\bigg)^2 dxds
  + C_3\int_0^t\|v_1-v_2\|_{L^2(\Omega)}^2 ds.
\end{align*}
\end{lemma}

\begin{proof}
We choose the test function $\phi=v_1-v_2$ in the weak formulation of the difference $v_1-v_2$ (see Lemma \ref{lem.weak2}):
\begin{align*}
  & \frac12\int_\Omega(v_1-v_2)^2(t)dx 
  - \frac12\int_\Omega(v_1-v_2)^2(0)dx = I_7 + I_8 + I_9,
  \quad\mbox{where} \\
  & I_7 = -\int_0^t\int_\Omega A_{22}(u_1,v_1)|\na(v_1-v_2)|^2 dxds, \\
  & I_8 = -\int_0^t\int_\Omega\big(A_{22}(u_1,v_1)-A_{22}(u_2,v_2)\big)
  \na v_2\cdot\na(v_1-v_2) dxds, \\
  & I_9 = \int_0^t\int_\Omega\big(R_2(v_1,v_1)-R_2(u_2,v_2)\big)
  (v_1-v_2)dxds.
\end{align*}
We infer from Assumption (A3) and the boundedness of $u_1$ and $v_1$ that $A_{22}$ is bounded from below by a positive constant $c_1$. This yields
\begin{align*}
  I_7 \le -c_1\int_0^t\int_\Omega|\na(v_1-v_2)|^2 dxds.
\end{align*}
Similarly as for the term $I_{51}$ in the proof of Lemma \ref{lem.ineq2}, we estimate
\begin{align*}
  I_8 &\le C\int_0^t\int_\Omega\big(|u_1^{\alpha/2+1}-u_2^{\alpha/2+1}|
  + |v_1-v_2|\big)|\na(v_1-v_2)|dxds \\
  &\le C\int_0^t\int_\Omega\Big(\sqrt{(u_1^{\alpha+1}-u_2^{\alpha+1})
  (u_1-u_2)} + |v_1-v_2|\Big)|\na(v_1-v_2)|dxds \\
  &\le \frac{c_1}{4}\int_0^t\int_\Omega|\na(v_1-v_2)|^2 dxds
  + C\int_0^t\int_\Omega(u_1^{\alpha+1}-u_2^{\alpha+1})
  (u_1-u_2)dxds \\
  &\phantom{xx}+ C\int_0^t\int_\Omega(v_1-v_2)^2 dxds.
\end{align*}
The remaining term $I_9$ is estimated similarly as $I_6$ in the proof of Lemma \ref{lem.ineq2}. Indeed, we write $I_9 = I_{91}+I_{92}+I_{93}$, where
\begin{align*}
  I_{91} &= \int_0^t\int_\Omega(u_1-u_2)q_2(v_1)(v_1-v_2)dxds, \\
  I_{92} &= \int_0^t\int_\Omega u_2(q_2(v_1)-q_2(v_2))(v_1-v_2)dxds, \\
  I_{93} &= \int_0^t\int_\Omega\big(R_2^0(u_1,v_1)-R_2^0(u_2,v_2)\big)
  (v_1-v_2)dxds.
\end{align*}
Inserting equation \eqref{2.psi}, integrating by parts, and using Young's inequality as for the estimate of $I_{61}$ in the proof of Lemma \ref{lem.ineq2}, we obtain
\begin{align*}
  I_{91} &= \int_0^t\int_\Omega q_2'(v_1)(v_1-v_2)\na(\psi_1-\psi_2)
  \cdot\na v_1dxds \\
  &\phantom{xx}+ \int_0^t\int_\Omega q_2(v_1)\na(\psi_1-\psi_2)
  \cdot\na(v_1-v_2)dxds \\
  &\phantom{xx}+ \int_0^t\int_\Omega q_2(v_1)(v_1-v_2)dx
  \bigg(\fint_\Omega(u_1-u_2)dx\bigg)ds \\
  &\le C\int_0^t\int_\Omega|\na(\psi_1-\psi_2)|^2 dxds
  + C\int_0^t\int_\Omega(v_1-v_2)^2 dxds \\
  &\phantom{xx}+ \frac{c_1}{4}\int_0^t\int_\Omega|\na(v_1-v_2)|^2 dxds
  + C\int_0^t\bigg(\int_\Omega(u_1-u_2)dx\bigg)^2 ds.
\end{align*}
By similar arguments as above, we have
\begin{align*}
  I_{92}+I_{93} \le C\int_0^t\int_\Omega(u_1^{\alpha+1}-u_2^{\alpha+1})
  (u_1-u_2)dxds + C\int_0^t\int_\Omega(v_1-v_2)^2 dxds. 
\end{align*}
We conclude that
\begin{align*}
  I_9 &\le \frac{c_1}{4}\int_0^t\int_\Omega|\na(v_1-v_2)|^2 dxds
  + C\int_0^t\int_\Omega|\na(\psi_1-\psi_2)|^2 dxds \\
  &\phantom{xx}+ C\int_0^t\int_\Omega(v_1-v_2)^2 dxds
  + C\int_0^t\int_\Omega(u_1^{\alpha+1}-u_2^{\alpha+1})(u_1-u_2)dxds \\
  &\phantom{xx}+ C\int_0^t\bigg(\int_\Omega(u_1-u_2)dx\bigg)^2ds. 
\end{align*}
Combining the estimate for $I_7$, $I_8$, and $I_9$ finishes the proof.
\end{proof}

We proceed with the proof of Theorem \ref{thm.cont}. We need to absorb the terms involving $\na(v_1-v_2)$ and $(u_1^{\alpha+1}-u_2^{\alpha+1}) (u_1-u_2)$. For this, we multiply the inequality of Lemma \ref{lem.ineq2} by the factor $K>0$ and add the inequalities of Lemmas \ref{lem.ineq1} and \ref{lem.ineq3}. Then
\begin{align*}
  \bigg(\int_\Omega&(u_1-u_2)(t)dx\bigg)^2
  + K\|\na(\psi_1-\psi_2)(t)\|_{L^2(\Omega)}^2
  + \|(v_1-v_2)(t)\|_{L^2(\Omega)}^2 \\
  &\phantom{xx}+ (c_1-\delta K)
  \int_0^t\|\na(v_1-v_2)\|_{L^2(\Omega)}^2 ds \\
  &\phantom{xx}+ (K c_0 - C_1 - C_3)\int_0^t\int_\Omega 
  (u_1^{\alpha+1}-u_2^{\alpha+1}) (u_1-u_2)dxds \\
  &\le 4\bigg(\int_\Omega(u_1^0-u_2^0)^2 dx\bigg)^2
  + \|v_1^0-v_2^0\|_{L^2(\Omega)}^2
  + K\|\na(\psi_1-\psi_2)(0)\|_{L^2(\Omega)}^2 \\
  &\phantom{xx}+ (C_1+KC_2+C_3)\int_0^t\bigg\{
  \bigg(\int_\Omega(u_1-u_2)dx\bigg)^2
  + \delta^{-1}\|\na(\psi_1-\psi_2)\|_{L^2(\Omega)}^2 \\
  &\phantom{xx}+ \|v_1-v_2\|_{L^2(\Omega)}^2\bigg\}ds.
\end{align*}
The fourth and fifth terms on the left-hand side are both nonnegative if we choose $K\ge(C_1+C_3)/c_0$ and $\delta<\min\{1,c_1/K\}$.

We claim that 
\begin{align}\label{3.psi0}
  \int_\Omega|\na(\psi_1-\psi_2)(0)|^2 dx
  \le C\int_\Omega(u_1^0-u_2^0)^2 dx.
\end{align}
Indeed, using the test function $(\psi_1-\psi_2)(0)$ in the weak formulation of \eqref{2.psi} for $i=1,2$ at time $t=0$, we find that
\begin{align*}
  \int_\Omega |\na(\psi_1-\psi_2)(0)|^2 dx
  &= \int_\Omega(u_1^0-u_2^0)(\psi_1-\psi_2)(0)dx \\
  &\phantom{xx}- \fint_\Omega(u_1^0-u_2^0)dx
  \int_\Omega(\psi_1-\psi_2)(0)dx.
\end{align*}
Then a standard application of the Poincar\'e and Young inequalities leads to
\begin{align*}
  \int_\Omega |\na(\psi_1-\psi_2)(0)|^2 dx
  \le C\int_\Omega(u_1^0-u_2^0)^2dx
  + \frac12\int_\Omega|\na(\psi_1-\psi_2)(0)|^2 dx,
\end{align*}
which shows claim \eqref{3.psi0}. The proof of inequality \eqref{1.cont} follows from 
\begin{align*}
  \|u_1-u_2\|_{H^1(\Omega)'}^2 
  \le C\|\na(\psi_1-\psi_2)\|_{L^2(\Omega)}^2
  + C\bigg(\int_\Omega(u_1-u_2)dx\bigg)^2
\end{align*}
for some constant $C>0$. This inequality follows from definition of the $H^1(\Omega)'$ norm and definition \eqref{2.psi} of $\psi_i$,
\begin{align*}
  \|u_1-u_2\|_{H^1(\Omega)'}
  &= \sup_{\|\phi\|_{H^1(\Omega)}=1}|\langle u_1-u_2,\phi\rangle| \\
  &= \sup_{\|\phi\|_{H^1(\Omega)}=1}\bigg|
  \int_\Omega\bigg(\na(\psi_1-\psi_2)\cdot\na\phi 
  + \fint_\Omega(u_1-u_2)dx\phi\bigg)dx\bigg| \\
  &\le \|\na(\psi_1-\psi_2)\|_{L^2(\Omega)} 
  + C\bigg|\int_\Omega(u_1-u_2)dx\bigg|.
\end{align*}
Finally, the uniqueness result is a direct consequence of inequality \eqref{1.cont}.


\end{document}